\documentclass[runningheads]{llncs}
\usepackage{graphicx}
\usepackage{amsmath}
\usepackage{amsfonts}
\usepackage{bm}
\usepackage{mathtools}
\usepackage{stmaryrd}
\usepackage{tikz,pgfplots}
\pgfplotsset{compat=newest}
\usetikzlibrary{positioning,calc,patterns,math}
\usepackage{siunitx}
\usepackage{url}

\newcommand{\bs}[1]{\boldsymbol{#1}}
\newcommand{\ue}{u_\mathrm{e}}
\newcommand{\Vm}{V_\mathrm{m}}

\newcommand{\vX}{\mathbf{X}}

\newcommand{\vx}{\mathbf{x}}

\newcommand{\vn}{\mathbf{n}}
\newcommand{\tA}{\mathbf{a}}
\newcommand{\tG}{\mathbf{G}}

\newcommand{\tGi}{\tG_\mathrm{i}}
\newcommand{\tGe}{\tG_\mathrm{e}}
\newcommand{\tGb}{\tG_\mathrm{b}}
\newcommand{\tGm}{\tG_\mathrm{m}}

\newcommand{\R}{\mathbb{R}}
\newcommand{\dd}{\mathrm{d}}
\newcommand{\EE}{\mathbb{E}}
\newcommand{\OmegaH}{\Omega_\mathrm{H}}
\newcommand{\OmegaT}{\Omega_\mathrm{T}}

\newcommand{\sil}{\sigma_{\mathrm{i},\mathrm{l}}}
\newcommand{\sit}{\sigma_{\mathrm{i},\mathrm{t}}}
\newcommand{\sel}{\sigma_{\mathrm{e},\mathrm{l}}}
\newcommand{\set}{\sigma_{\mathrm{e},\mathrm{t}}}

\newcommand{\Vrest}{V_\mathrm{rest}}

\newcommand{\Vdepol}{V_\mathrm{dep}}

\DeclareMathOperator{\Cor}{Cor}
\DeclareMathOperator{\Var}{Var}

\begin{document}
\title{Fast and Accurate Uncertainty Quantification
for the ECG with Random Electrodes Location}
%
%%%%%%%%%%%%%%%%%%%%%%%%%%%%%%%%%%%%%%%%%%%%%%%%%%%%%%%%%%%%%%%%%%%%%%%%%%%%%%%%
\titlerunning{UQ in ECG with lead field approach}
\author{Michael Multerer\orcidID{0000-0003-0170-0239} \and Simone Pezzuto\orcidID{0000-0002-7432-0424}}
\authorrunning{Multerer \and Pezzuto}
%%%%%%%%%%%%%%%%%%%%%%%%%%%%%%%%%%%%%%%%%%%%%%%%%%%%%%%%%%%%%%%%%%%%%%%%%%%%%%%%
\institute{Center for Computational Medicine in Cardiology,
Institute of Computational Science,
Universit\`a della Svizzera italiana,
via G.~Buffi 13, 6900 Lugano, Switzerland
\email{\{michael.multerer,simone.pezzuto\}@usi.ch}}
\maketitle
%
%%%%%%%%%%%%%%%%%%%%%%%%%%%%%%%%%%%%%%%%%%%%%%%%%%%%%%%%%%%%%%%%%%%%%%%%%%%%%%%%
\begin{abstract}
The standard electrocardiogram (ECG) is a point-wise evaluation of the body
potential at certain given locations.  These locations are subject to uncertainty
and may vary from patient to patient or even for a single patient. In this work,
we estimate the uncertainty in the ECG induced by uncertain electrode positions
when the ECG is derived from the forward bidomain model.
In order to avoid the high computational cost associated to the solution of the
bidomain model in the entire torso, we propose a low-rank approach to solve the
uncertainty quantification (UQ) problem.
More precisely, we exploit the sparsity of the ECG and the lead field theory
to translate it into a set of deterministic, time-independent problems,
whose solution is
eventually used to evaluate expectation and covariance of the ECG. We assess
the approach with numerical experiments in a simple geometry.
%%%%%%%%%%%%%%%%%%%%%%%%%%%%%%%%%%%%%%%%%%%%%%%%%%%%%%%%%%%%%%%%%%%%%%%%%%%%%%%%
\keywords{Random Electrodes Location
\and Uncertainty Quantification
\and Lead Field
\and Electrophysiology
\and Forward Bidomain Model}
\end{abstract}
%%%%%%%%%%%%%%%%%%%%%%%%%%%%%%%%%%%%%%%%%%%%%%%%%%%%%%%%%%%%%%%%%%%%%%%%%%%%%%%%

\section{Introduction}
%%%%%%%%%%%%%%%%%%%%%%%%%%%%%%%%%%%%%%%%%%%%%%%%%%%%%%%%%%%%%%%%%%%%%%%%%%%%%%%%
The standard ECG is a routinely acquired recording of the torso electric
potential~\cite{plonsey95}. 
%It consists in 12 recordings of the electric potential
%on the chest.  Each recording or \emph{lead} is the potential difference between
%two electrodes or between an electrode and a reference potential.
%For instance,
%lead I is the potential difference between left and right arm, whereas lead V1
%is obtained from the difference of an electrode on the chest (called again V1)
%and the Wilson Central Terminal, that is the average of limb electrode potentials.
It provides valuable information on the electric activity
of the heart and, when combined with imaging data of the anatomy,
it can be used for non-invasive personalization of sophisticated
patient-specific models~\cite{GiffardRoisin2019,Pezzuto2021}.
In these inverse ECG models, the ECG is rarely computed from the state-of-the-art
bidomain model~\cite{colli2014}, otherwise the computational cost would be prohibitive. Commonly,
the bidomain model is replaced by a ``decoupled'' version, called
forward bidomain~\cite{potse2006} or pseudo-bidomain~\cite{bishop2011,neic2017} model,
in which the transmembrane potential in the heart is computed independently
from the extracellular potential in the torso.
The resulting model still compares favourably to the coupled bidomain model and,
more importantly, the ECG can be evaluated
very efficiently and \emph{exactly} by employing the lead field theory~\cite{pezzuto2017,potse2018}.

Obviously, when dealing with real data, as in patient-specific modeling,
model parameters are subject to unavoidable uncertainty. This uncertainty
should be accounted for
in the forward and inverse ECG model~\cite{clayton2020}.
%Uncertainty in the parameters turns the simulated ECG into a stochastic
%variable and renders the forward model a possibly expensive
%uncertainty quantification problem, due to the
%high non-linearity in the stochastic variables and the
%intrinsic cost in a single model evaluation.
Several sources of
uncertainty may be considered, e.g., related to the segmentation
process of the anatomy~\cite{corrado2020},
the electric conductivities~\cite{aboulaich2016},
or the fiber distribution~\cite{Quaglino2018}.
Particularly relevant in the context of inverse ECG modeling
is the uncertainty in the electrodes' locations, which has shown
to yield sensible morphological changes in the precordial signals
even with a displacement as low as \SI{2}{\cm}~\cite{kania2014}.

The present work focuses on the problem of estimating the expectation and
the covariance of the surface ECG, if electrodes' locations are subject to uncertainty
and the ECG is simulated with the forward bidomain model.
In principle, given the torso potential, the statistical moments are readily available
with little additional cost, as the solution of the UQ problem amounts
to a simple integration
over the torso domain.  In spite of its simplicity, the computational
cost of this approach grows linearly with the number of time steps and
the number of evaluations of the forward model. Moreover, it
relies on the full torso potential, despite the fact that the electrodes'
locations may be very localized.
We propose a computationally very efficient methodology to
solve the UQ problem \emph{without} the need of solving the full forward problem.
Our method is still based on the lead field theory and it is an
\emph{exact} representation of the true ECG. Specifically,
it exploits a low-rank approach
to decouple the correlation problem into a small set of elliptic
problems for different right hand sides~\cite{HL13}.
Remarkably, the overall computational cost is drastically reduced
and comparable to the solution of a few
elliptic problems, \emph{independently} of the number of time steps and forward
evaluations.
%Instead of solving the space-time problem, we recast it into a limited number of simpler time-independent problems by a low-rank approximation of the covariance function.
%The proposed methodology is even more beneficial in a multi-query context, e.g., when multiple instances of the forward problem, with different transmembrane potential but same torso model, are to be computed.

This paper is organized as follows: in Sec.~\ref{sec:meth},
we review the forward bidomain model for the ECG,
the lead field approach and describe our method.
In Sec.~\ref{sec:exp}, we validate the approach on a simple geometry.
We conclude in Sec.~\ref{sec:conc} with a brief discussion and outlook.

%Making the assumption, that the electrode positions are very localized, we
%employ a low-rank approach to 
%decouple the tensor product lead field problem for the correlation into a
%set of lead field problems for different right hand sides,
%see \cite{%Har14,
%HL13}. If the low-rank approach is not applicable, an
%alternative would be to resort to hierarchical
%matrix techniques, see \cite{%DHP15,
%DHP17}.
%The presented approach is heavily inspired by the lead field theory~\cite{potse2018}. More specifically, we focus on estimating the expectation and the correlation function of the ECG when this is resulting from a bidomain simulation in the whole torso.

\section{Methods}
\label{sec:meth}

%%%%%%%%%%%%%%%%%%%%%%%%%%%%%%%%%%%%%%%%%%%%%%%%%%%%%%%%%%%%%%%%%%%%%%%%%%%%%%%%
\subsection{The forward bidomain model}
The electric potential $u_0(\vx,t)$ in the torso $\OmegaT\subset\R^d$,
and consequently the ECG, can be modelled from the transmembrane potential
$\Vm(\vx,t)$ in the active myocardium $\OmegaH\subset\R^d$,
with the time-dependent \emph{forward bidomain model}~\cite{potse2006},
which reads as follows:
%%%%%%%%%%%%%%%%%%%%%%%%%%%%%%%%%%%%%%%%%%%%%%%%%%%%%%%%%%%%%%%%%%%%%%%%%%%%%%%%
\begin{equation}\label{eq:bidomain}
\begin{cases}
%\beta\Bigl(\Cm\partial_t\Vm + \fion(\Vm,\vz)\Bigr) = \nabla\cdot(\tGm\nabla\Vm),
%& \mbox{in $\OmegaH\times[0,\infty)$}, \\
-\nabla\cdot\bigl((\tGi+\tGe)\nabla\ue(\vx,t)\bigr)  = \nabla\cdot(\tGi\nabla\Vm(\vx,t)),
& \mbox{in $\OmegaH\times[0,\infty)$}, \\
-\nabla\cdot(\tG_0\nabla u_0(\vx,t)) = 0, & \mbox{in $\OmegaT\times[0,\infty)$}, \\
-\tG_0\nabla u_0(\vx,t)\cdot\vn = 0, & \mbox{on $\Sigma\times[0,\infty)$}, \\
\phantom{-}\ue(\vx,t) = u_0(\vx,t), & \mbox{on $\Gamma\times[0,\infty)$}, \\
-\tGb\nabla\ue(\vx,t)\cdot\vn + \tG_0\nabla u_0(\vx,t)\cdot\vn = \tGi\nabla\Vm(\vx,t)\cdot\vn, &
\mbox{on $\Gamma\times[0,\infty)$}.
\end{cases}
\end{equation}
%%%%%%%%%%%%%%%%%%%%%%%%%%%%%%%%%%%%%%%%%%%%%%%%%%%%%%%%%%%%%%%%%%%%%%%%%%%%%%%%

Herein, $\Gamma = \bar\Omega_\mathrm{H}\cap\bar\Omega_\mathrm{T}$ is the
heart-torso interface, $\Sigma = \partial\OmegaT\setminus\Gamma$ is the body
surface, $\ue(\vx,t)$ is the extra-cellular potential in the heart,
$\tGi$ and $\tGe$ are respectively intra- and extra-cellular conductivity
of the heart, $\tG_0$ is the torso
conductivity, and $\vn$ is the outward normal for both $\Gamma$ and $\Sigma$.
For the sake of simplicity in the notation, we define
\[
\tG \coloneqq
\begin{cases}
\tGi+\tGe &\text{in }\OmegaH, \\
\tG_0 &\text{in }\OmegaT,
\end{cases}
\quad
u(\vx,t) \coloneqq
\begin{cases}
\ue(\vx,t) & \text{in }\OmegaH, \\
u_0(\vx,t) & \text{in }\OmegaT,
\end{cases}
\]
%%%%%%%%%%%%%%%%%%%%%%%%%%%%%%%%%%%%%%%%%%%%%%%%%%%%%%%%%%%%%%%%%%%%%%%%%%%%%%%%
and assume, without loss of generality, that
$u(\cdot,t)\in\mathrm{H}^1(\Omega)$, where $\Omega = \OmegaH \cup \OmegaT$.
In this case, the variational formulation for Eq.~\eqref{eq:bidomain}
can be written according to 
\begin{equation}\label{eq:varform}
\begin{aligned}
&\text{For every $t\in\R$, find $u(\cdot,t)\in\mathrm{H}^1(\Omega)$
such that}\\
&\quad\int_{\Omega} \tG\nabla u(\vx,t) \cdot\nabla v \:\dd\vx
= - \int_{\OmegaH}
\tGi \nabla\Vm(\vx,t) \cdot\nabla v \:\dd\vx
\end{aligned}
\end{equation}
%%%%%%%%%%%%%%%%%%%%%%%%%%%%%%%%%%%%%%%%%%%%%%%%%%%%%%%%%%%%%%%%%%%%%%%%%%%%%%%%
for all $v\in \mathrm{H}^1(\Omega)$. The well-posedness of the problem follows from
standard application of the Riesz Theorem~\cite{evans2010},
given that $\OmegaH, \OmegaT$ are Lipschitz domains and
$\Vm(\cdot,t)\in\mathrm{H}^1(\OmegaH)$.  We remark that the formulation in Eq.~\eqref{eq:varform} is equivalent to Eq.~\eqref{eq:bidomain} when the
restriction of the solution $u|_{\Omega_i}$ belongs to $\mathrm{H}^2(\Omega_i)$,
\(i\in\{\text{H},\text{T}\}\), see e.g.~\cite{ammari2016,chen1998} for a more 
comprehensive
treatment of interface problems.

The ECG is a set of so-called \emph{leads}, typically 12 in the standard ECG.
Each lead reads as follows:
\begin{equation}\label{eq:ecg}
V(t,{\bs\xi}_1,\ldots,{\bs\xi}_L) = \sum_{\ell=1}^L a_\ell u({\bs\xi}_\ell,t),
\end{equation}
where ${\bs\Xi}\coloneqq \{{\bs\xi}_\ell \}_{\ell=1}^L$ is the set of electrodes
and $\tA = [a_1,\ldots,a_L]^\top$ is a zero-sum vector of weights defining the
lead. For instance, a limb lead is the potential difference of 2 electrodes,
whereas a precordial lead involves 4 electrodes (3 are used to build the
Wilson Central Terminal, that is the reference potential).
It is worth noting that Eq.~\eqref{eq:ecg} is valid only
if $u(\cdot,t)\in\mathcal{C}^0(\Sigma)$, which is not true for
$u(\cdot,t)\in\mathrm{H}^1(\Omega)$ and $d\ge 2$. For a rigorous discussion, see~\cite{colli2014}.

%%%%%%%%%%%%%%%%%%%%%%%%%%%%%%%%%%%%%%%%%%%%%%%%%%%%%%%%%%%%%%%%%%%%%%%%%%%%%%%%
In this work, we are interested in computing statistics of
$V\big(t,{\bs\Xi}(\omega)\big)$ when
the electrode positions
${\bs\Xi(\omega)\coloneqq\{{\bs\xi}_\ell(\omega)\}_{\ell=1}^L}$ are not known exactly.
Here, we
denote by ${\bs\xi}_\ell(\omega)$ the random variable associated to the $\ell$-th
electrode and assume that the joint distribution is given by the density
$\rho(\vX) = \rho(\vx_1,\ldots,\vx_L)$ with respect to the surface measure
$\dd\sigma_{\vX} = \dd\sigma_{\vx_1}\cdots\dd\sigma_{\vx_L}$
on $\Sigma^L$.
According to the definition in Eq.~\eqref{eq:ecg}, the lead $V(t,{\bs\Xi})$
is a random field as well, with expectation and correlation respectively
reading as follows:
\begin{align}
\EE[V](t) &= \int_{\Sigma^M} V(t,\vX) \rho(\vX)\dd\sigma_{\vX},
\label{eq:avg} \\
\Cor[V](t,s) &= \int_{\Sigma^M} V(t,\vX) V(s,\vX)\rho(\vX)
\dd\sigma_{\vX}. \label{eq:cov}
\end{align}

In summary, the UQ problem for the random electrodes locations
consists in solving the forward bidomain model Eq.~\eqref{eq:bidomain}
for $u(\vx,t)$, given $\Vm(\vx,t)$, for every time $t$, and then computing
the integrals in Eq.~\eqref{eq:avg} and Eq.~\eqref{eq:cov}.

\subsection{Lead field formulation of the UQ problem}
%%%%%%%%%%%%%%%%%%%%%%%%%%%%%%%%%%%%%%%%%%%%%%%%%%%%%%%%%%%%%%%%%%%%%%%%%%%%%%%%
Clearly, in general it is not convenient to compute the ECG from Eq.~\eqref{eq:bidomain},
because the ECG is only a very sparse evaluation of $u(\vx,t)$.  Moreover, in
a patient-specific or personalization context, the ECG needs to be simulated
several times with different instances of $\Vm(\vx,t)$, with no changes in the
left hand side of Eq.~\eqref{eq:bidomain}. A better approach is based on
Green's functions, also known as \emph{lead fields} in the
electrocardiographic literature~\cite{potse2018}. In fact, it is possible to
show that $V(t,{\bs\Xi})$ has the following representation~\cite{colli2014}:
\begin{equation}\label{eq:ecglead}
V(t,{\bs\Xi}) = \int_{\OmegaH} \tGi(\vx)\nabla\Vm(\vx,t) \cdot\nabla Z(\vx,{\bs\Xi})
\:\dd{\vx},
\end{equation}
where $Z(\vx,{\bs\Xi})$ is the weak solution of the elliptic problem:
\begin{equation}\label{eq:leadpb}
\begin{cases}
-\nabla\cdot\tG\nabla Z(\vx,{\bs\Xi}) = 0, & \text{in }\Omega, \\
-\tG\nabla Z(\vx,{\bs\Xi})\cdot \vn = \sum_{\ell=1}^L a_\ell \delta_{{\bs\xi}_\ell},
& \text{on }\Sigma,
\end{cases}
\end{equation}
%%%%%%%%%%%%%%%%%%%%%%%%%%%%%%%%%%%%%%%%%%%%%%%%%%%%%%%%%%%%%%%%%%%%%%%%%%%%%%%%
where $\delta_{{\bs\xi}_\ell}$ is the $(d-1)$-dimensional Dirac delta centered
at ${\bs\xi}_\ell(\omega)$. Therefore, given that all measurement locations are fixed,
Eq.~\eqref{eq:leadpb} is only solved once, at the cost of a single time step of
Eq.~\eqref{eq:bidomain}, and then used to compute $V(t,{\bs\Xi})$ for any choice
of $\Vm(\vx,t)$.
%To note, Eq.~\eqref{eq:ecglead} is often adopted with
%approximated choices of $Z({\bs\Xi})$, e.g., assuming an isotropic and
%homogeneous torso via the boundary element method, or even simpler the
%fundamental solution of the Laplacian in $\R^3$. Finally, we remark that
%Eq.~\eqref{eq:leadpb} is pure Neumann problem. Hence, the solution is only
%determined up to a constant. To fix this constant, we shall focus on solutions
%that vanish on average.

%\subsection{Expected ECG.}
Here, we exploit Eq.~\eqref{eq:ecglead} to compute the the expectation
and correlation of $V$, according to Eq.~\eqref{eq:avg} and \eqref{eq:cov}.
Substituting Eq.~\eqref{eq:ecglead} into Eq.~\eqref{eq:avg}, we obtain by
the linearity of the expectation that
\begin{equation}\label{eq:avg_lead}
\EE[V](t) = \int_{\OmegaH} \tGi(\vx)\nabla\Vm(\vx,t) \cdot\nabla \EE[Z](\vx) \:\dd{\vx}.
\end{equation}
Again by linearity, the equation for the expected lead field $\EE[Z]$ follows
from Eq.~\eqref{eq:leadpb} and reads as follows:
\begin{equation}\label{eq:leadpb_avg}
\begin{cases}
-\nabla\cdot\tG\nabla\EE[Z](\vx) = 0, & \text{in }\Omega, \\
-\tG\nabla \EE[Z](\vx)\cdot \vn = \sum_{\ell=1}^M a_\ell \rho_\ell(\vx),
& \text{on }\Sigma,
\end{cases}
\end{equation}
where $\rho_\ell$ is the marginal distribution of $\rho$ with respect to ${\bs\xi}_\ell$, that is
\begin{equation} \label{eq:margin1}
\rho_\ell(\vx_\ell) \coloneqq
\int_{\Sigma^{L-1}} \rho(\vX)\:\dd\sigma_{\vx_1}
\cdots \dd\sigma_{\vx_{\ell-1}}\dd\sigma_{\vx_{\ell+1}}\cdots \dd\sigma_{\vx_L}.
\end{equation}
To show this, we observe that:
\[
\EE\biggl[\sum_{\ell=1}^L a_\ell \delta_{{\bs\xi}_\ell}\biggr] =
\sum_{i=1}^L a_\ell \int_{\Sigma^L} \delta_{\vx_\ell} \rho(\vX)
\;\dd\sigma_{\vX} \\
= \sum_{\ell=1}^L a_\ell \rho_\ell.
\]
%%%%%%%%%%%%%%%%%%%%%%%%%%%%%%%%%%%%%%%%%%%%%%%%%%%%%%%%%%%%%%%%%%%%%%%%%%%%%%%%
Therefore, the cost of computing the average ECG is equivalent to that for
solving for the point-wise ECG, i.e., one solution of the elliptic problem
in Eq.~\eqref{eq:leadpb_avg}.
We observe that both Eq.~\eqref{eq:leadpb} and Eq.~\eqref{eq:leadpb_avg}
are well-posed, since
the right hand side has zero average over $\Sigma$ in both cases.
In particular, for every $\omega$, $Z(\vx,{\bs\Xi}(\omega))$ 
and $\EE[Z](\vx)$ are only defined up to a constant.

%Moreover, the averaged Eq.~\eqref{eq:leadpb_avg} is more regular than
%the original problem with no singularity at the electrodes location. We remark
%however that the presence of the singularity is irrelevant for the computation
%of the ECG, since Eq.~\eqref{eq:ecglead} considers the restriction of $Z$ on
%$\OmegaH$, clearly not in contact with the chest. However, the same approach
%may be applied to simulate intracardiac or intramural signals as well,
%in which case special care is required.\medskip

%\subsection{Correlation of the ECG.}
%%%%%%%%%%%%%%%%%%%%%%%%%%%%%%%%%%%%%%%%%%%%%%%%%%%%%%%%%%%%%%%%%%%%%%%%%%%%%%%%
The natural continuation of the above argument yields the correlation for
the ECG according to
\[
\Cor[V](t,s) 
= \int_{\Sigma^2} (\tGi\nabla\otimes\tGi\nabla) \Vm(\vx,t)\Vm(\vx',s)
: (\nabla\otimes\nabla) \Cor[Z]\: \dd\sigma_{\vx}\dd\sigma_{\vx'},
\]
where the tensor product is $[\mathbf{u}\otimes\mathbf{v}]_{ij} = u_i(\vx)
v_j(\vx')$ and the inner product between tensors is $\mathbf{A}:\mathbf{B}
= \sum_{ij} [A]_{ij} [B]_{ij}$.  The problem for the correlation $\Cor[Z]$,
obtained as above from Eq.~\eqref{eq:leadpb}, reads as follows:
\begin{equation}\label{eq:leadpb_cor}
\begin{cases}
(\nabla\cdot\tG\nabla\otimes \nabla\cdot\tG\nabla) \Cor[Z] = 0, & \text{in }
\Omega\times\Omega, \\
(\vn\cdot\tG\nabla \otimes \nabla\cdot\tG\nabla) \Cor[Z] = 0, &
\text{on } \Sigma\times\Omega, \\
(\nabla\cdot\tG\nabla \otimes \vn\cdot\tG\nabla) \Cor[Z] = 0, &
\text{on } \Omega\times\Sigma, \\
(\vn\cdot\tG\nabla \otimes \vn\cdot\tG\nabla) \Cor[Z] = R, &
\text{on }\Sigma\times\Sigma,
\end{cases}
\end{equation}
where the correlation $R(\vx,\vx')$ of the Neumann data in
Eq.~\eqref{eq:leadpb} is
\begin{equation}
\begin{split}
R(\vx,\vx') &= \Cor\biggl[\sum_{\ell=1}^L a_\ell \delta_{{\bs\xi}_\ell},
\sum_{\ell'=1}^L a_{\ell'} \delta_{{\bs\xi}_{\ell'}} \biggr] \\
&= \sum_{\ell=1}^L a_\ell^2 \rho_\ell(\vx)\delta_{\vx}(\vx')
+ \sum_{\ell\neq\ell'}^L a_\ell a_{\ell'} \rho_{\ell,\ell'}(\vx,\vx'),
\end{split}
\end{equation}
with $\rho_{\ell,\ell'}(\vx,\vx')$ being the marginal distribution
of $\rho$ with respect to $(\bs\xi_\ell,\bs\xi_\ell')$ and defined
as follows:
\begin{equation} \label{eq:margin2}
\rho_{\ell,\ell'}(\vx_\ell,\vx_{\ell'}) \coloneqq
\int_{\Sigma^{L-2}} \rho(\vX)\:\dd\sigma_{\vx_1}
\cdots \dd\sigma_{\vx_{\ell-1}}\dd\sigma_{\vx_{\ell+1}}
\cdots \dd\sigma_{\vx_{\ell'-1}}\dd\sigma_{\vx_{\ell'+1}}\cdots
\dd\sigma_{\vx_L}.
\end{equation}
We observe that, when $\bs\xi_\ell$ and $\bs\xi_\ell'$ are independent,
$\rho_{\ell,\ell'}(\vx_\ell,\vx_{\ell'})$ factorizes into the product
of the marginals $\rho_\ell(\vx_\ell)$ and $\rho_{\ell'}(\vx_{\ell'})$.

%%%%%%%%%%%%%%%%%%%%%%%%%%%%%%%%%%%%%%%%%%%%%%%%%%%%%%%%%%%%%%%%%%%%%%%%%%%%%%%%
As the computation of $\Cor[Z]$ requires the solution of a tensor product
boundary value problem, it is computationally rather expensive. 
In what follows, we will exploit the particular structure of $R(\vx,\vx')$
to significantly reduce the computational cost and implementation effort.
We remark that also Eq.~\eqref{eq:leadpb_cor} is well-posed because, by construction,
$R(\vx,\vx')$ is such that $\langle R,1\otimes v\rangle_{\Sigma^2}=
\langle R,v\otimes 1\rangle_{\Sigma^2}
= 0$ for all \(v\in\mathrm{H}^{1/2}(\Sigma)\),
where $\langle \cdot, \cdot \rangle_{\Sigma^2}$
is the duality pairing in \(\mathrm{L}^2(\Sigma^2)\).

\subsection{Numerical Discretization}
%%%%%%%%%%%%%%%%%%%%%%%%%%%%%%%%%%%%%%%%%%%%%%%%%%%%%%%%%%%%%%%%%%%%%%%%%%%%%%%%
The variational formulation of the averaged lead field
problem Eq.~\eqref{eq:avg_lead} resembles Eq.~\eqref{eq:varform} with a
different right hand side. With $Y=\mathrm{H}^1(\Omega)$, the problem is:
\begin{align*}
&\text{Find $\EE[Z]\in Y$ such that}\\[-0.5em]
&\qquad\int_{\Omega} \tG\nabla \EE[Z] \cdot\nabla v \:\dd\vx = 
\int_{\Sigma} \sum_{\ell=1}^L a_\ell \rho_\ell v\:\dd\vx,
\quad\text{for all $v\in Y$}.  
\end{align*}
The Galerkin approximation in the space $Y_h\subset Y$,
with $Y_h = \operatorname{span}\{ \phi_k\}_{k=1}^{N_h}$,
reads as follows:
\begin{equation}\label{eq:Kzg}
\mathbf{K}\mathbf{z} = \mathbf{g},
\end{equation}
where $\mathbf{z}$ is the solution vector, that is $\EE[Z]\approx Z_h = \sum_k 
[\mathbf{z}]_k \phi_k$ and
\begin{align}
[\mathbf{K}]_{k\ell} &= \int_{\Omega} \tG\nabla \phi_\ell\cdot\nabla \phi_k
\:\dd\vx, \label{eq:Kassemble} \\
[\mathbf{g}]_{k} &= \int_{\Sigma}
\sum_{\ell=1}^M a_\ell\rho_\ell(\vx)\phi_k(\vx)\:\dd\vx. \label{eq:gassemble}
\end{align}
%%%%%%%%%%%%%%%%%%%%%%%%%%%%%%%%%%%%%%%%%%%%%%%%%%%%%%%%%%%%%%%%%%%%%%%%%%%%%%%%
For the correlation in Eq.~\eqref{eq:leadpb_cor}, the variational
formulation is as follows:
\begin{align*}
&\text{Find $\Cor[Z]\in Y\otimes Y$ such that}\\[-0.25em]
&\qquad\int_{\Omega^2} (\tG\nabla\otimes\tG\nabla) \Cor[Z] :
(\nabla\otimes\nabla) v \:\dd\vx\dd\vx' = 
\int_{\Sigma^2} Rv \:\dd\vx\dd\vx'
\end{align*}
for all $v\in Y\otimes Y$.  The corresponding Galerkin formulation on
$Y_h \times Y_h$ is:
\begin{equation}\label{eq:GalerkinCor}
(\mathbf{K}\otimes\mathbf{K})\mathbf{Z} = \mathbf{R},
\end{equation}
where $\Cor[Z]\approx \sum_{k,\ell} [\mathbf{Z}]_{k\ell}
\phi_k\otimes\phi_\ell$
and
\begin{equation}\label{eq:Rassemble}
\begin{split}
[\mathbf{R}]_{pq} &= \int_{\Sigma^2} R\phi_p\phi_q \dd\vx\dd\vx' \\
&= \sum_{\ell=1}^L a_\ell^2 \int_\Sigma \rho_\ell\phi_p \phi_q \dd\vx
+ \sum_{\ell\neq\ell'}^L a_\ell a_{\ell'}
\int_{\Sigma^2} \rho_{\ell,\ell'}\phi_p \phi_q \dd\vx\dd\vx'.
%\left( \int_\Sigma \rho_\ell \phi_p \dd\vx \right)\left( \int_\Sigma
%\rho_{\ell'} \phi_q \dd\vx \right).
\end{split}
\end{equation}

%%%%%%%%%%%%%%%%%%%%%%%%%%%%%%%%%%%%%%%%%%%%%%%%%%%%%%%%%%%%%%%%%%%%%%%%%%%%%%%%
As the number of degrees of freedom for the correlation problem is \(N_h^2\),
it may easily become computationally prohibitive. However, assuming that the
marginal densities \(\rho_\ell\), \(\ell=1,\ldots,L\) are strongly localized, the
right hand side in \eqref{eq:GalerkinCor} may be represented by a low-rank
approximation according to
\[
\mathbf{R}\approx \sum_{k=1}^K \bm{r}_k \otimes \bm{r}_k,
\quad\bm{r}_k\in\mathbb{R}^{N_h},
\]
with \(K\ll N_h\).
In this case, we also expect a low-rank solution, that is
\[
\mathbf{Z} \approx \sum_{k=1}^K \bm{\zeta}_k \otimes \bm{\zeta}_k,
\quad\bs{\zeta}_k\in\mathbb{R}^{N_h},
\]
with $K \ll N_h$. 
Then, due to the tensor product structure of \eqref{eq:GalerkinCor}, there
simply holds
\begin{equation}\label{eq:Kzeta}
\mathbf{K} \bm{\zeta}_k = \bm{r}_k, \qquad k=1,\ldots,K.
\end{equation}
In practice, we compute the low-rank approximation by a diagonally
pivoted, truncated Cholesky decomposition, see \cite{HPS12}.

%%%%%%%%%%%%%%%%%%%%%%%%%%%%%%%%%%%%%%%%%%%%%%%%%%%%%%%%%%%%%%%%%%%%%%%%%%%%%%%%
Finally, for the computation of statistics of the ECG, we insert the computed
Galerkin approximations into Eq.~\eqref{eq:avg} and Eq.~\eqref{eq:cov} and
obtain
\begin{align}
\EE[V](t)    &\approx \mathbf{V}(t)\cdot \mathbf{z}, \label{eq:EEV} \\
\Cor[V](t,s) &\approx \sum_{k,m=1}^K\big(\mathbf{V}(t)
\cdot \bm{\zeta}_k\big)\big(\mathbf{V}(s)\cdot \bm{\zeta}_m\big), \label{eq:CorV}
\end{align}
where
\[
[\mathbf{V}(t)]_j = \int_{\OmegaH} \tGi\nabla\Vm(t)\cdot \nabla\phi_j\:\dd\vx.
\]

The computational cost for the proposed approach is dominated
by the solution of $K+1$
systems (one for the expectation and $K$ for the correlation)
of the form of Eq.~\eqref{eq:Kzg}. It is therefore
independent on the number of time steps $N_t$ or the choice
of $\Vm$, oppositely to the
solution of the forward bidomain model in Eq.~\eqref{eq:bidomain},
which requires $K\cdot N_t$ solutions for each choice of $\Vm$.

We summarize below the proposed procedure to evaluate expectation and correlation of a single lead
of the ECG, defined with $L \ge 2$ coefficients $\mathbf{a}$
as in Eq.~\eqref{eq:ecglead} and with random electrodes locations
$\{ \bs\xi_\ell \}_{\ell=1}^L$
with density $\rho(\vX)$. We assume as above that $\Vm(\vx,t)$ is given
and computed elsewhere.
\begin{enumerate}
\item Compute $\rho_\ell(\vx)$ with Eq.~\eqref{eq:margin1} and assemble $\mathbf{g}$ with Eq.~\eqref{eq:gassemble};
\item Assemble $\mathbf{K}$ with Eq.~\eqref{eq:Kassemble} and solve Eq.~\eqref{eq:Kzg} to find $\mathbf{z}$;
\item Compute $\EE[V](t)$ from $\mathbf{z}$ and $\Vm$ with Eq.~\eqref{eq:EEV};
\item Compute $\rho_{\ell,\ell'}(\vx,\vx')$ with Eq.~\eqref{eq:margin2} and assemble $\mathbf{R}$ with Eq.~\eqref{eq:Rassemble};
\item Compute the low-rank Cholesky decomposition $\{ \bs{r}_k \}_{k=1}^K$ of $\mathrm{R}$;
\item For each $k=1,\ldots,K$, solve Eq.~\eqref{eq:Kzeta} for $\bs\zeta_k$;
\item Compute $\Cor[V](t,s)$ from $\{ \bs{\zeta}_k \}_{k=1}^K$ and $\Vm$ with Eq.~\eqref{eq:CorV}.
\end{enumerate}

\section{Numerical Assessment}
\label{sec:exp}

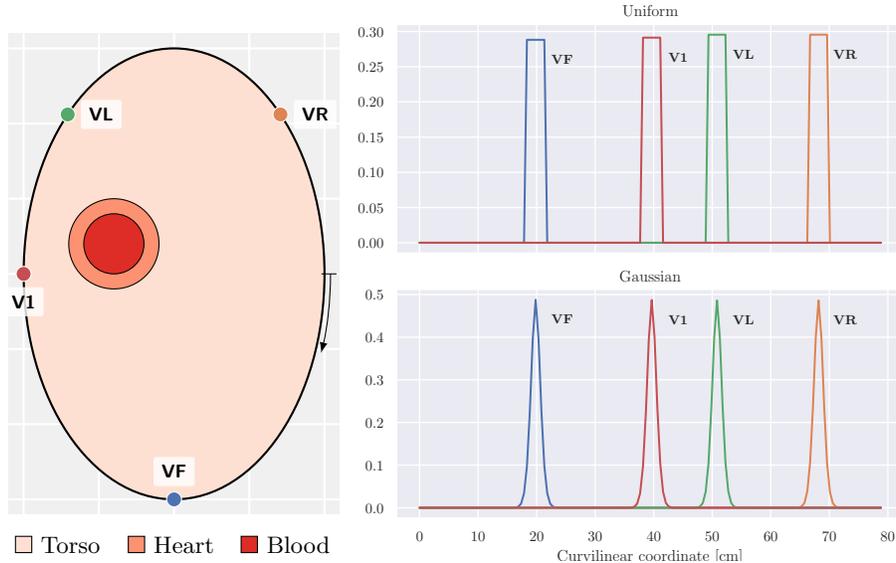
\begin{figure}[tb]
    \centering
    \begin{tikzpicture}
    \def\myunit{2mm}
    \definecolor{colheart}{HTML}{fc9272}
    \definecolor{coltorso}{HTML}{fee0d2}
    \definecolor{colblood}{HTML}{de2d26}
    \definecolor{colVF}{HTML}{4c72b0}
    \definecolor{colVR}{HTML}{dd8452}
    \definecolor{colVL}{HTML}{55a868}
    \definecolor{colV1}{HTML}{c44e52}

    \definecolor{colgrid}{HTML}{f0f0f0}
    \tikzstyle{elec}=[circle,inner sep=2pt,draw=white,thin]
    \tikzstyle{heart}=[colheart,even odd rule,draw=black]
    \tikzstyle{torso}=[coltorso,draw=black,thick]
    \tikzstyle{tbox}=[rounded corners=1pt,fill=white,opacity=0.8,text opacity=1.0,
                      outer sep=2pt,font=\bfseries\sffamily\scriptsize]
    \tikzset{
    partial ellipse/.style args={#1:#2:#3}{
        insert path={+ (#1:#3) arc (#1:#2:#3)}
    }}
    % grid
    \fill[colgrid] (-11*\myunit,-16*\myunit) rectangle (11*\myunit,16*\myunit);
    \draw [step=5*\myunit,white,thick] (-11*\myunit,-16*\myunit) grid (11*\myunit,16*\myunit);
    % torso
    \fill[torso]  (0,0) ellipse [x radius = 10*\myunit, y radius = 15*\myunit];
    % heart
    \fill[colblood] (-4*\myunit,2*\myunit) circle[radius = 2*\myunit];
    \fill[heart] (-4*\myunit,2*\myunit) circle[radius = 3*\myunit] circle[radius = 2*\myunit];
    % electrodes
    \node[elec,fill=colVL] (VL) at ({10*\myunit*cos(deg(3*pi/4))},{15*\myunit*sin(deg(3*pi/4))}) {};
    \node[elec,fill=colVR] (VR) at ({10*\myunit*cos(deg(1*pi/4))},{15*\myunit*sin(deg(1*pi/4))}) {};
    \node[elec,fill=colVF] (VF) at ({10*\myunit*cos(deg(6*pi/4))},{15*\myunit*sin(deg(6*pi/4))}) {};
    \node[elec,fill=colV1] (V1) at ({10*\myunit*cos(deg(4*pi/4))},{15*\myunit*sin(deg(4*pi/4))}) {};
    % labels
    \node[tbox,right] at (VL.east) {VL};
    \node[tbox,right] at (VR.east) {VR};
    \node[tbox,above] at (VF.north) {VF};
    \node[tbox,below] at (V1.south) {V1};
    % arc-length
    \draw[thin] (9.8*\myunit,0) -- +(1.0*\myunit,0);
    \draw[thin,latex-] (0,0) [partial ellipse=-20:0:{10.4*\myunit} and {15.4*\myunit}];
    % legend
    \path (-10*\myunit,-18*\myunit) node[draw=black,fill=coltorso,label=right:{Torso}] (b1) {}
        +(7.5*\myunit,0) node[draw=black,fill=colheart,label=right:{Heart}] (b2) {}
       ++(15*\myunit,0) node[draw=black,fill=colblood,label=right:{Blood}] (b3) {};
    
    % distributions
    % \node[anchor=north west] (dist1) at (15*\myunit,16*\myunit)
    % {\includegraphics[scale=0.5,trim=50 30 40 30,clip]{rho_dist}};
    % \node[anchor=south] at (dist1.north) {Uniform distribution};

    % \node[anchor=north west] (dist2) at ($(dist1.south west)+(0,-5mm)$)
    % {\includegraphics[scale=0.5,trim=50 0 40 30,clip]{rho_smooth}};
    % \node[anchor=south] at (dist2.north) {Gaussian-like distribution};
    
    % \node[anchor=center]
    % at (dist2.south) {\scriptsize Curvilinear coordinate [cm]};
    
    \end{tikzpicture}\hfill
    \includegraphics[scale=0.5]{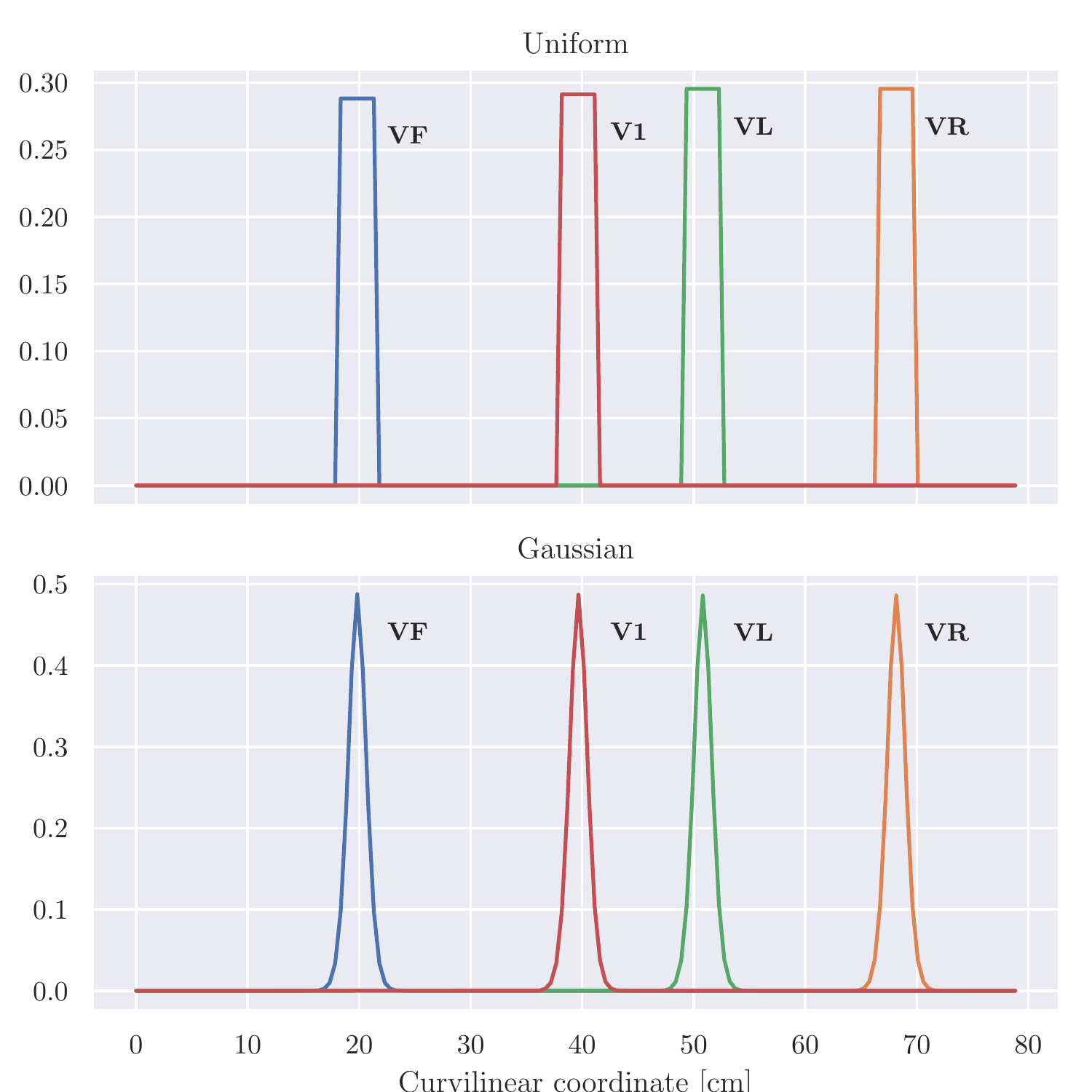}
    
    \caption{Geometrical configuration for the numerical test. On the left,
    the domain is represented with electrodes locations on the boundary and
    tissue properties. On the right, the probability density function for both
    uniform and Gaussian-like cases is
    reported. The $x$-axis is the curvilinear coordinate, and colors of
    the curves refer to the electrodes on the left.}
    \label{fig:geom2d}
\end{figure}

%%%%%%%%%%%%%%%%%%%%%%%%%%%%%%%%%%%%%%%%%%%%%%%%%%%%%%%%%%%%%%%%%%%%%%%%%%%%%%%%
We tested the proposed approach on a idealized heart-torso geometry in
2-D, as depicted in Fig.~\ref{fig:geom2d}.  The anatomy consists of an
ellipsoidal torso with major semi-axis of $T_y = \SI{15}{\cm}$, vertically oriented,
and minor axis of $T_x = \SI{10}{\cm}$.  The heart was an annulus centered at
$\vx_\mathrm{h}=(\SI{-4}{\cm}, \SI{2}{\cm})$ with respect to the center of the torso,
and with inner (endocardium) and outer (epicardium) radius respectively
equal to \SI{2}{\cm} and \SI{3}{\cm}.  The domain was split into 3 distinct
regions, namely blood pool, myocardium and torso (see Fig.~\ref{fig:geom2d}).

For this test, we considered an ECG with 2 leads obtained from 4 random electrodes
$\bs\xi_\ell(\omega)$,
$\ell = \{ \mathrm{VL},\mathrm{VR},\mathrm{VF},\mathrm{V1} \}$,
see Fig.~\ref{fig:geom2d}. The leads, II and V1, were
\begin{align*}
V_\mathrm{II}(t,\bs\Xi)  &= u(\bs\xi_\mathrm{VF},t) - u(\bs\xi_\mathrm{VL},t), \\
V_\mathrm{V1}(t,\bs\Xi) &= u(\bs\xi_\mathrm{V1},t)
- \frac{1}{3}\Bigl( u(\bs\xi_\mathrm{VL},t) + u(\bs\xi_\mathrm{VR},t) + u(\bs\xi_\mathrm{VF},t) \Bigr),
\end{align*}
respectively corresponding to $\mathbf{a}_\mathrm{II} = (-1,0,1,0)$ and
$\mathbf{a}_\mathrm{V1} = (-\frac{1}{3},-\frac{1}{3},-\frac{1}{3},1)$.
The average position $\vx_\ell$ of the electrodes
was conveniently defined using the formula $\vx_\ell = [T_x \cos(\theta_\ell), T_y\sin(\theta_\ell)]^\top$
with $\theta_\mathrm{VL} = \frac{3}{4}\pi$, $\theta_\mathrm{VR} = \frac{1}{4}\pi$,
$\theta_\mathrm{VF} = \frac{3}{2}\pi$ and $\theta_\mathrm{V1} = \pi$.

In all tests, we evaluated the deterministic ECG,
computed from Eq.~\eqref{eq:ecglead}, the average ECG when electrodes were randomly
located, and the variance from the formula
$\Var[V](t) = \Cor[V](t,t) - \big(\EE[V](t)\big)^2$.

The random electrode locations, independent from each other,
were either uniformly distributed or
with a Gaussian-like distribution, both defined on the outer boundary $\Sigma$ of
the torso (the ``chest''), see Fig.~\ref{fig:geom2d}.  In the case of the
uniform distribution, the marginal density $\rho_\ell$ for each electrode was the characteristic function of the set $\Sigma_\ell = \{ \vx\in\Sigma: d(\vx_\ell,\vx) \le r_\ell \}$,
that is the $r_\ell$-neighborhood of $\vx_\ell$ with respect to the geodesic distance $d$ on the curve $\Sigma$.
We also considered a Gaussian-like distribution computed, after a normalization,
by solving the heat
equation on the boundary curve $\Sigma$ with diffusion defined along the
arc-length, initial datum $\delta(\vx-\vx_\ell)$ on $\Sigma$,
and solved for a total time $T_\ell$, see Fig.~\ref{fig:geom2d}.
We selected $r_\ell = \SI{1.5}{\cm}$ and $T_\ell = \sqrt{3}/3 r_\ell$,
so that both distributions have the same variance.

%%%%%%%%%%%%%%%%%%%%%%%%%%%%%%%%%%%%%%%%%%%%%%%%%%%%%%%%%%%%%%%%%%%%%%%%%%%%%%%%
%%%%%%%%%%%%%%%%%%%%%%%%%%%%%%%%%%%%%%%%%%%%%%%%%%%%%%%%%%%%%%%%%%%%%%%%%%%%%%%%

For convenience, we report the full expression $R(\vx,\vx')$
for lead II, obtained by assuming that $\bs\xi_\mathrm{VF}$
and $\bs\xi_\mathrm{VL}$ were independent:
\[
R(\vx,\vx') = \bigl(\rho_\mathrm{VF}(\vx)+\rho_\mathrm{VL}(\vx)\bigr)\delta_{\vx}(\vx')
- \rho_\mathrm{VF}(\vx)\rho_\mathrm{VL}(\vx') - \rho_\mathrm{VL}(\vx)\rho_\mathrm{VF}(\vx').
\]
In particular, the assembly of the tensor $\mathbf{R}$ simplifies as well, with
no need of evaluating a double integral. In fact,
\[
[\mathbf{R}]_{pq} = 
\int_\Sigma \bigl(\rho_\mathrm{VF}+\rho_\mathrm{VL}\bigr) \bigl(\phi_p\phi_q - \phi_p - \phi_q\bigr) \dd\sigma_\vx.
\]

The electric conductivities were uniform and isotropic in the torso and
in the blood pool, and respectively set to \SI{2}{\milli\siemens\per\cm\square}
and \SI{6}{\milli\siemens\per\cm\square}.  The myocardium was assumed
transversely isotropic, with fibers $\mathbf{f}$ circularly oriented
and of unit length. Specifically:
\begin{align*}
\tGi &= \sit\mathbf{I} + (\sil-\sit)\mathbf{f}\otimes\mathbf{f}, \\
\tGe &= \set\mathbf{I} + (\sel-\set)\mathbf{f}\otimes\mathbf{f},
\end{align*}
and values set to
$\sil = \SI{3}{\milli\siemens\per\cm\square}$,
$\sit = \SI{0.3}{\milli\siemens\per\cm\square}$,
$\sel = \SI{3}{\milli\siemens\per\cm\square}$ and
$\set = \SI{1.2}{\milli\siemens\per\cm\square}$.

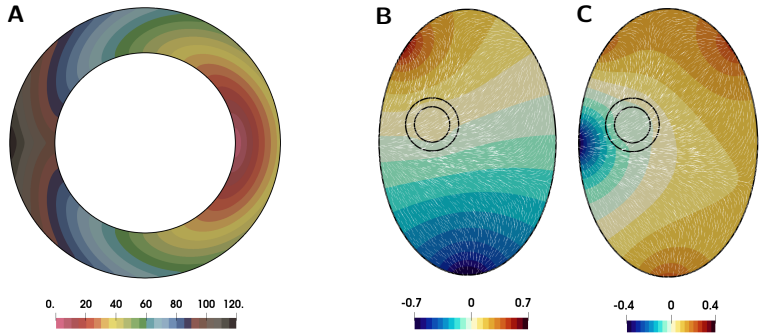
\begin{figure}[tb]
    \centering
    \begin{tikzpicture}
    \node (A) {\includegraphics[width=0.3\textwidth]{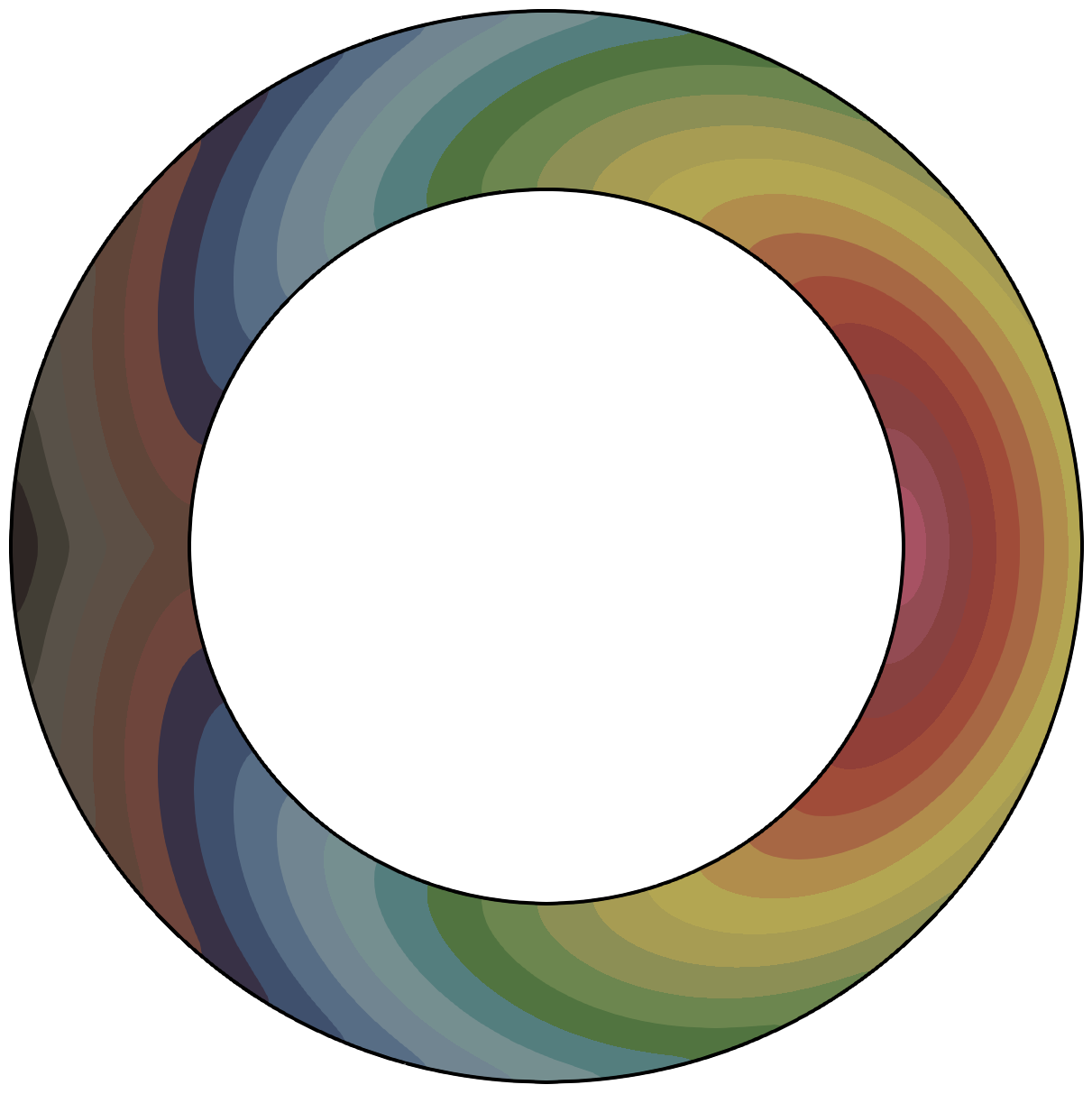}};
    \node[anchor=west] (B) at ($(A.east)+(1cm,0)$) {\includegraphics[width=0.2\textwidth]{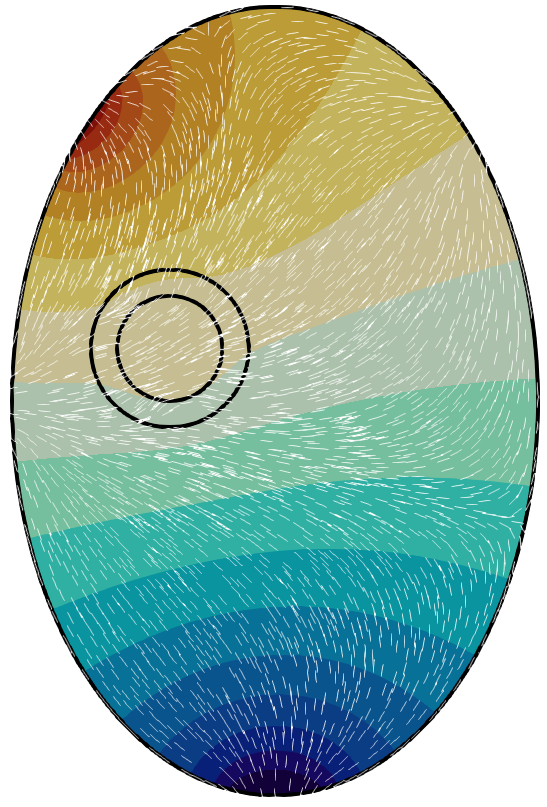}};
    \node[anchor=west] (C) at (B.east) {\includegraphics[width=0.2\textwidth]{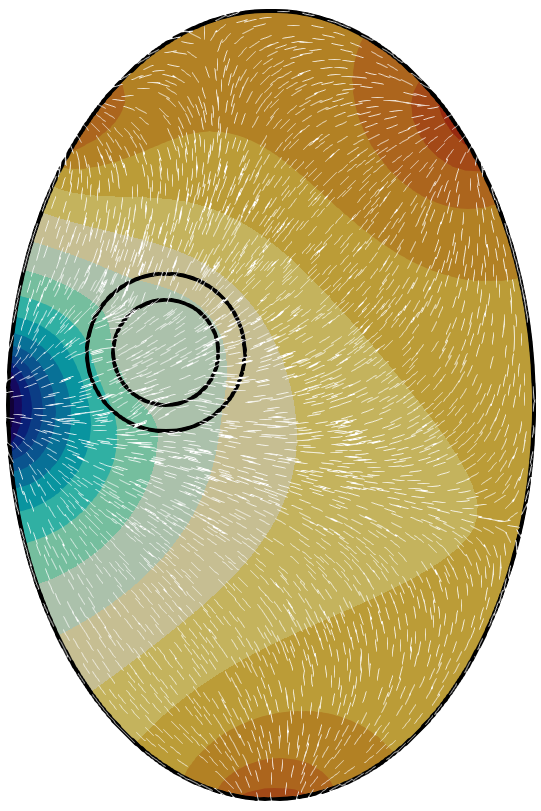}};
    % labels
    \node[anchor=north west,font=\bfseries\sffamily] at (A.north west) {A};
    \node[anchor=north west,font=\bfseries\sffamily] at (B.north west) {B};
    \node[anchor=north west,font=\bfseries\sffamily] at (C.north west) {C};
    % colorbars
    \node[anchor=north] (Al) at (A.south) {\includegraphics[width=0.22\textwidth]{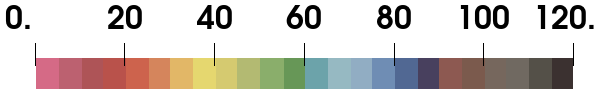}};
    \node[anchor=north] at (B.south) {\includegraphics[width=0.15\textwidth]{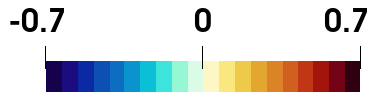}};
    \node[anchor=north] at (C.south) {\includegraphics[width=0.12\textwidth]{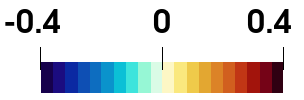}};
    \end{tikzpicture}
    \caption{On the left (panel A), the activation map computed with the eikonal solution, in \si{\ms}. On right, lead fields (in \si{\milli\volt}) for the average problem are reported for lead II (panel B) and lead V1 (panel C).}
    \label{fig:activ}
\end{figure}

For sake of simplicity,
the transmembrane potential $\Vm(\vx,t)$ was modelled by shifting a fixed action
potential template at given activation times $\tau(\vx)$, according to the
formula
\[
\Vm(\vx,t) = U\bigl(t-\tau(\vx)\bigr),
\qquad
U(s) = \Vrest + \frac{\Vdepol-\Vrest}{2}
\biggl(1 + \tanh\Bigl(\frac{s}{\varepsilon}\Bigr) \biggr).
\]
The activation map $\tau\colon \OmegaH\to\mathbb{R}$ was simulated with the
eikonal model
\[
\begin{cases}
\sqrt{ \mathbf{D}(\vx)\nabla\tau\cdot\nabla\tau } = 1, 
& \vx\in\OmegaH\setminus\{ \vx_\mathrm{s} \}, \\
\tau(\vx_\mathrm{s}) = 0. &
\end{cases}
\]
The conductivity tensor was set proportional to the monodomain conductivity
and such to yield a conduction velocity along the fibers of
\SI{65}{\cm\per\s}, that is,
\[
\mathbf{D} = \alpha\cdot\tGm,
\]
where $\tGm = \tGi(\tGi+\tGe)^{-1}\tGe$ and $\alpha \approx \SI{2.82e-3}{\cm\tothe{4}\square\ms\per\milli\siemens}$.
The other parameters were as follows:
$\Vrest = \SI{-85}{\milli\volt}$,
$\Vdepol = \SI{30}{\milli\volt}$,
$\varepsilon = \SI{0.4}{\ms}$,
and $\vx_\mathrm{s} = (-\SI{2}{\cm},\SI{2}{\cm})$.

%%%%%%%%%%%%%%%%%%%%%%%%%%%%%%%%%%%%%%%%%%%%%%%%%%%%%%%%%%%%%%%%%%%%%%%%%%%%%%%%
The computational domain was approximated by a triangular mesh $\mathcal{T}_h$
with median edge size of \SI{0.04}{\cm} in $\OmegaH$ and \SI{0.5}{\cm} in the rest of
the domain, thus totalling \num{27820} nodes and \num{55476} cells.
All quantities were represented by linear finite elements on $\mathcal{T}_h$.
The eikonal equation was solved with an anisotropic version of the heat distance
method~\cite{crane2013}, with $\Delta t = \SI{4}{\ms}$.
The implementation in \textsc{FEniCS} is publicly available\footnote{See \url{https://github.com/pezzus/fimh2021}.}
and complemented with additional tests and comparison to the monodomain
and bidomain models.

\begin{figure}[tb]
    \centering
    \includegraphics[width=\textwidth,trim=10 10 10 10,clip]{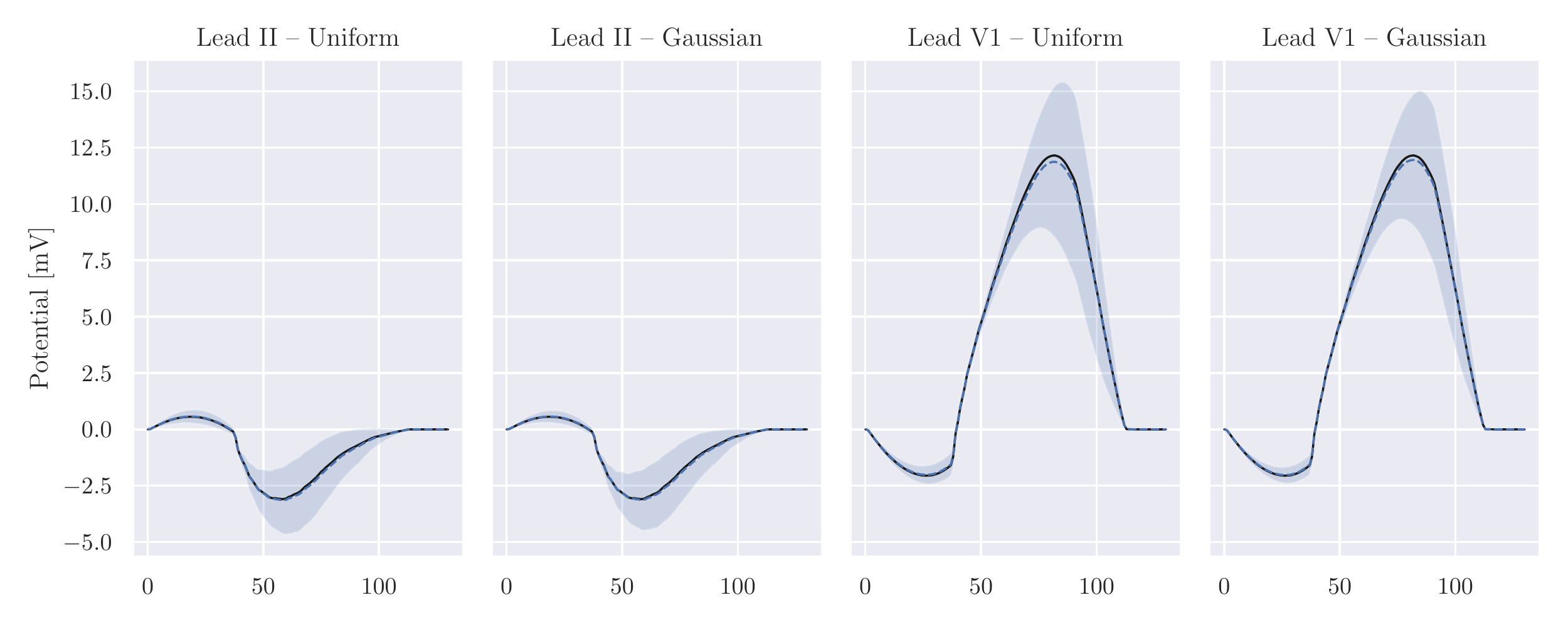}
    \caption{ECG results for all tests. In the plots, the dashed black curve
    is the deterministic ECG, the solid blue curve is the average
    ECG, and the shaded blue area corresponds to the \SI{95}{\percent} confidence
    interval, that is $\EE[V](t) \pm 1.96\sqrt{\Var[V](t)}$.}
    \label{fig:results}
\end{figure}

The activation map and the average lead fields for lead II and lead V1, as computed from Eq.~\eqref{eq:avg_lead}, are reported in Fig.~\ref{fig:activ}.
In the correlation problem, the low-rank representation counted 17 (resp.~33) modes for lead II with uniform (resp.~Gaussian) distribution of electrodes, and 33 (resp.~59) modes for lead V1. A lower number of modes for the uniform distribution was expected, as its support was compact and highly localized.
The resulting ECGs are reported in Fig.~\ref{fig:results}.  In both leads, the deterministic and average ECGs were very close, with an absolute error between \SI{0.083}{\milli\volt} (lead II) and
\SI{0.28}{\milli\volt} (lead V1). The uncertainty was significantly higher in the late part of the QRS-complex. Maximum standard deviation was as high as \SI{2.07}{\milli\volt} in lead V1 and \SI{0.83}{\milli\volt} in lead II.  In lead V1, the morphological variations were limited but the
maximum amplitude changed significantly. In lead II, morphological differences were present in the second half of the ECG.  No appreciable differences in ECGs were noted when comparing uniform and Gaussian-like distributions.  

\begin{figure}[tb]
    \centering
    \begin{tikzpicture}
    \node[inner sep=0pt,outer sep=0pt] (a) {\includegraphics[width=0.22\textwidth]{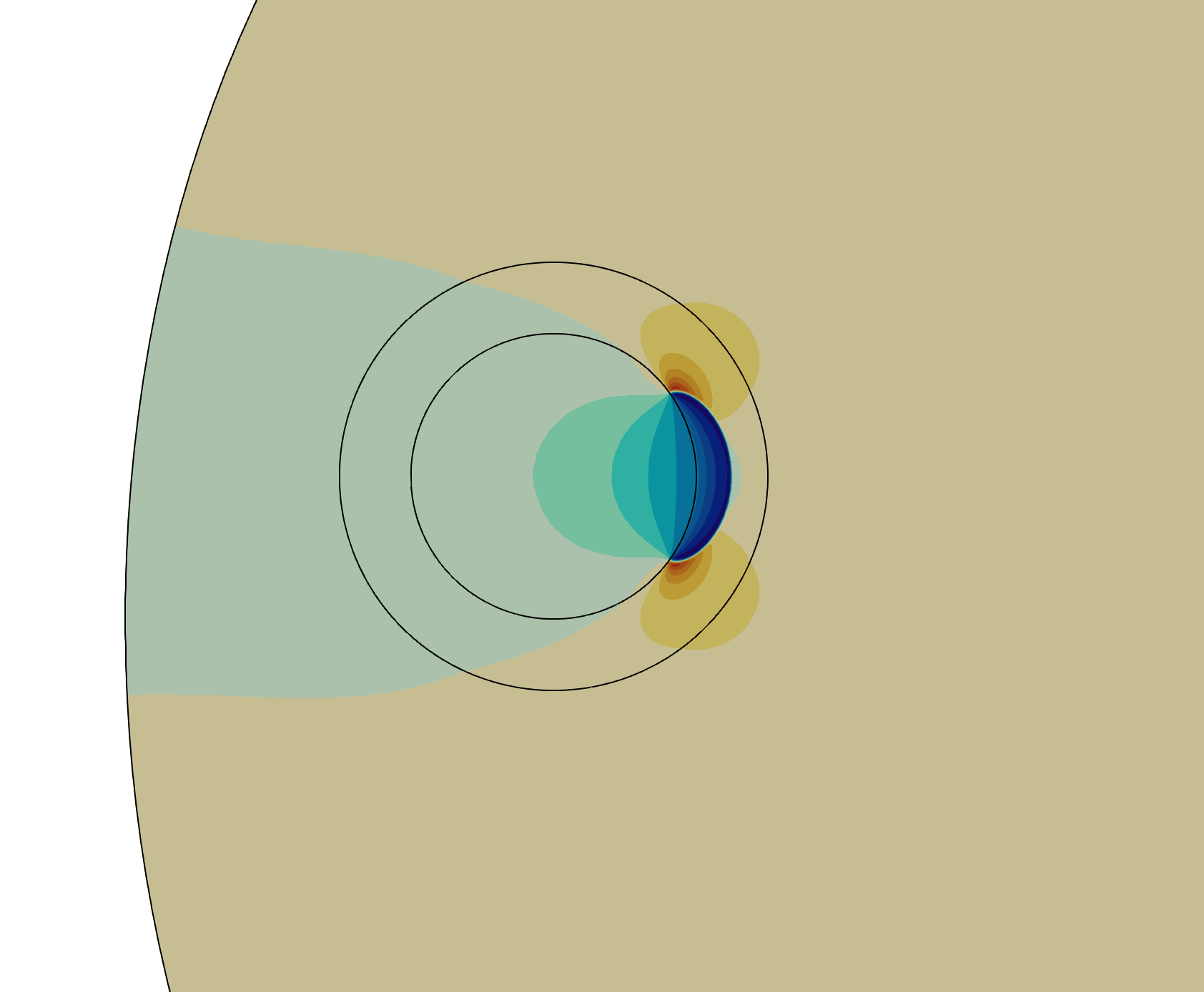}}
    node[anchor=south] at (a.north) {$t = \SI{20}{\ms}$};
    \end{tikzpicture}%
    %%%%
    \begin{tikzpicture}
    \node[inner sep=0pt,outer sep=0pt] (a) {\includegraphics[width=0.22\textwidth]{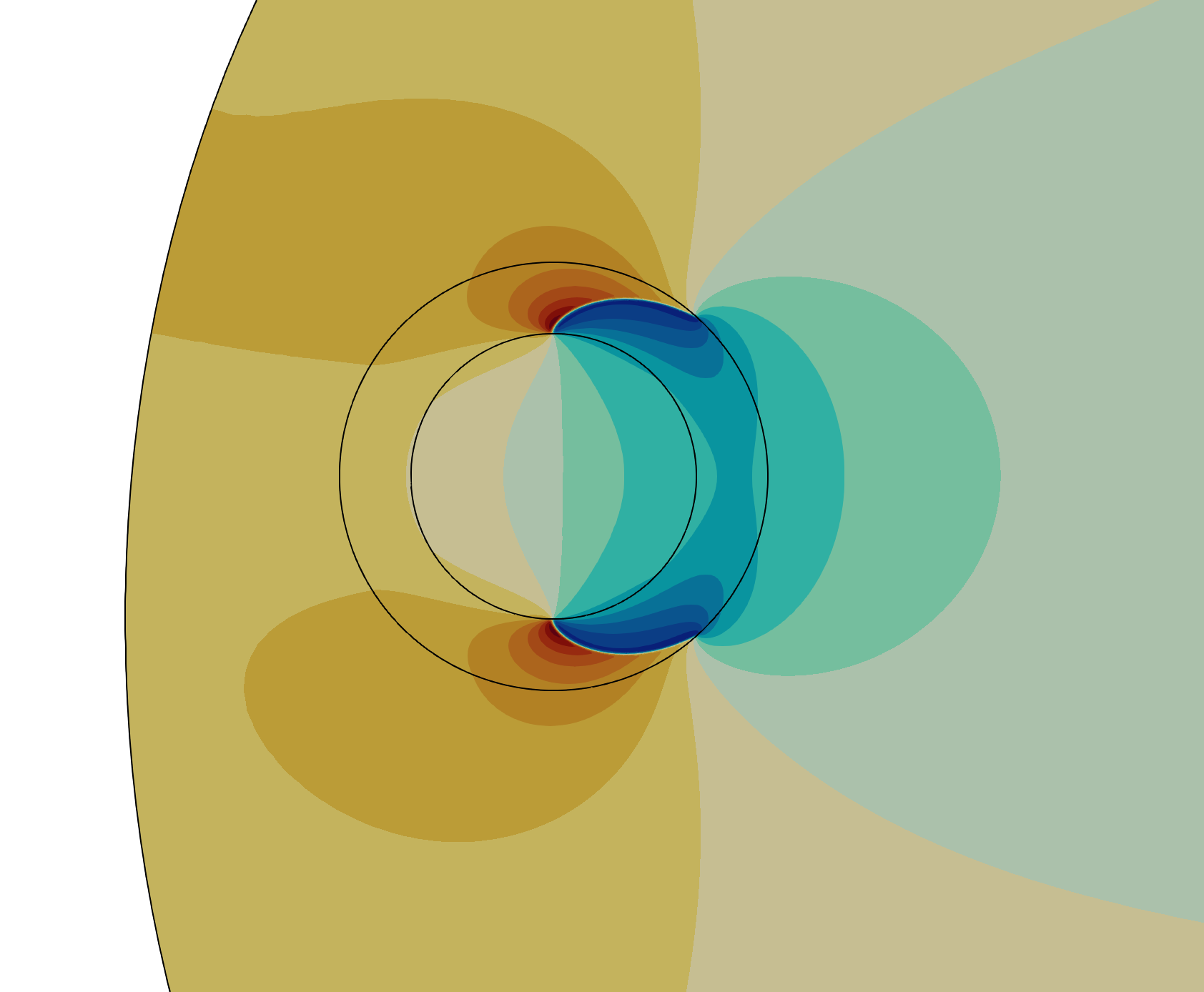}}
    node[anchor=south] at (a.north) {$t = \SI{50}{\ms}$};
    \end{tikzpicture}%
    %%%%
    \begin{tikzpicture}
    \node[inner sep=0pt,outer sep=0pt] (a) {\includegraphics[width=0.22\textwidth]{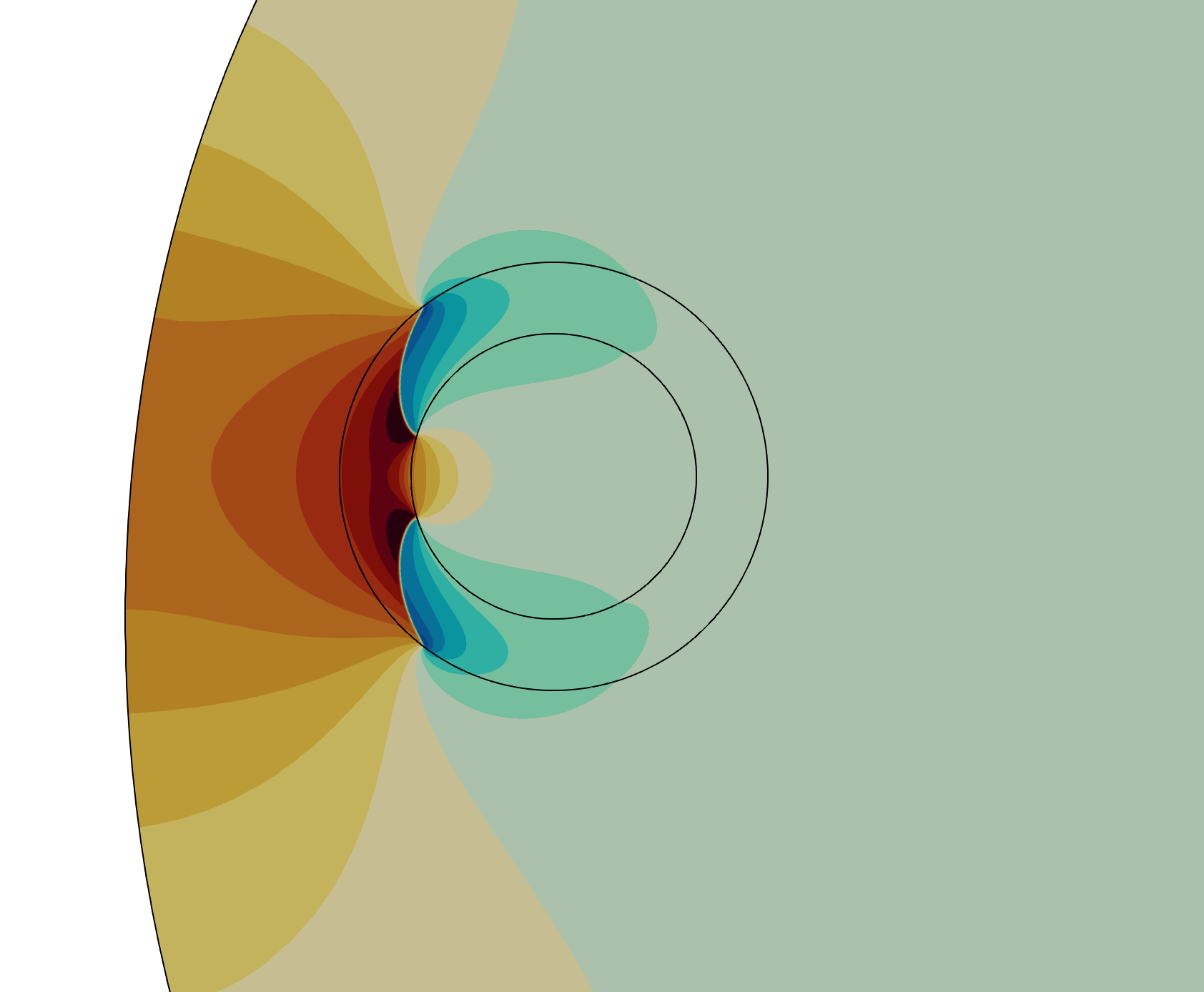}}
    node[anchor=south] at (a.north) {$t = \SI{90}{\ms}$};
    \end{tikzpicture}%
    %%%%
    \begin{tikzpicture}
    \node[inner sep=0pt,outer sep=0pt] (a) {\includegraphics[width=0.22\textwidth]{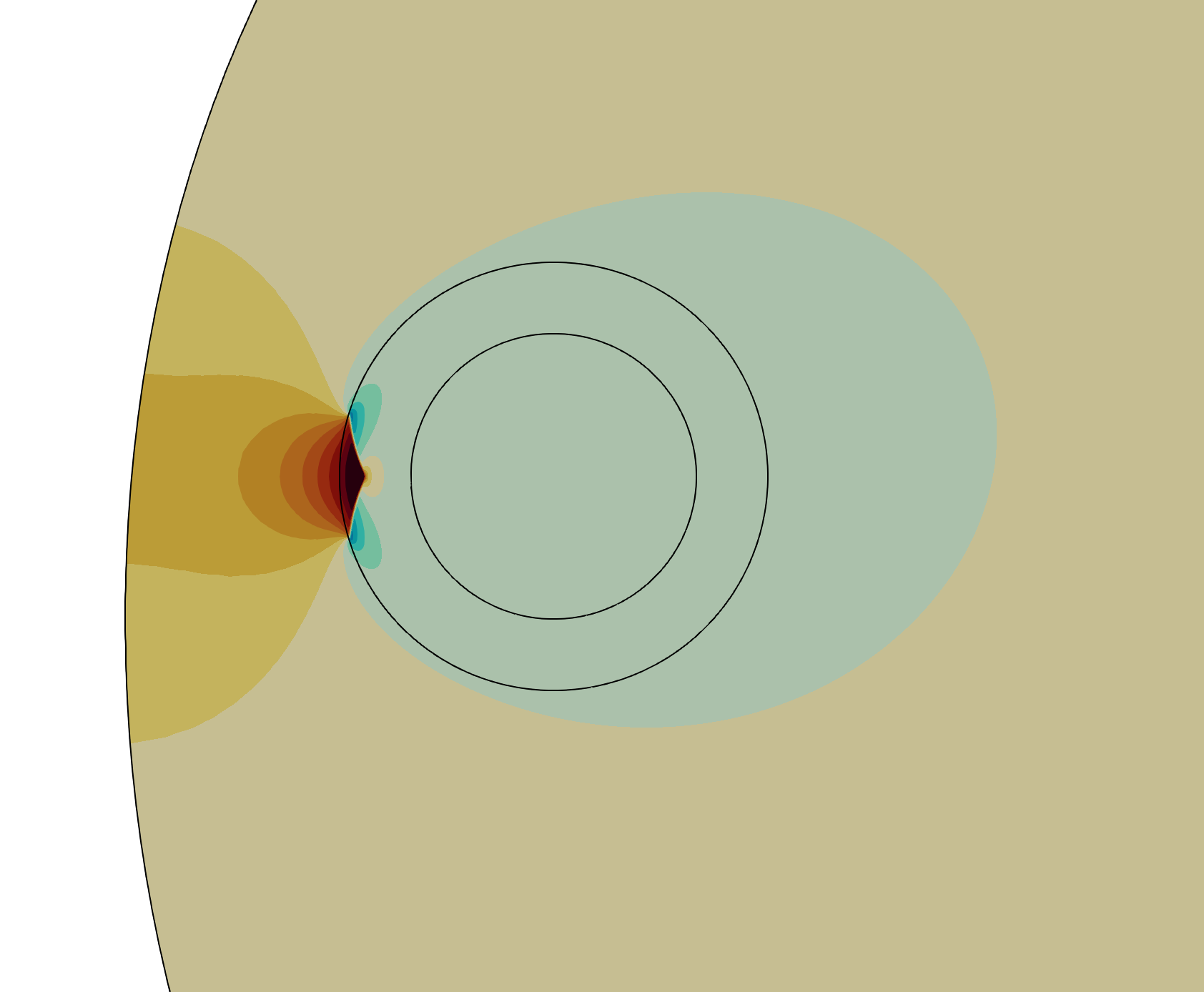}}
    node[anchor=south] at (a.north) {$t = \SI{110}{\ms}$};
    \end{tikzpicture}
    %%%%
    \includegraphics[width=0.06\textwidth]{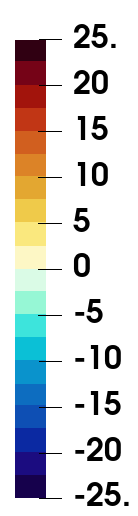}
    \caption{Excerpt of extracellular potential (in \si{\milli\volt})
    of the forward bidomain simulation.}
    \label{fig:bido}
\end{figure}

Finally, we compared the proposed method against
the solution of the forward bidomain model, see Fig.~\ref{fig:bido}.  Differences between ECGs computed from our approach were essentially matching those derived from the forward bidomain model, with an absolute error less than \SI{0.01}{\milli\volt} in all cases and quantities of interests. The total cost of the forward simulation was from 2 to 8-fold higher than the lead field approach.

\section{Discussion and Conclusions}
\label{sec:conc}

In this work, we have solved the problem of quantifying the uncertainty
in the ECG when uncertainty in the electrode positions is taken into account.
Our method recasts the
problem into a fully deterministic setting by using the lead field theory
and a low-rank approximation for the correlation.

The computational advantage
is significant over the standard forward simulation of the bidomain model.
In fact, the number of lead fields to be computed in the proposed approach,
for both the expectation and the correlation, does not depend on neither the 
transmembrane potential nor time, oppositely to the bidomain model.
The method is therefore
suitable to compute ECGs for long simulations, e.g., arrhythmic events,
and it is even more advantageous in the context of inverse ECG approaches.
Finally, the lead fields are smoother than the extracellular potential,
especially within the heart, where potential gradients are strong along the activation front.  A much coarser resolution may be employed for computing
the lead fields, with no significant loss in accuracy~\cite{potse2018}.

While formulated for the chest electrodes, the presented theory
also applies with minimal changes to assess the uncertainty of intracardiac
electrogram recordings, widely employed in clinical electrophysiological studies.
As a matter of fact, the formulation is flexible enough to address other relevant problems,
such as quantifying the uncertainty in the ECG due to, e.g., uncertain
transmembrane potential or torso-heart segmentation, hence leading
to more robust simulation results.

\bibliographystyle{splncs04}
\bibliography{randecg}
\end{document}